\documentclass[10pt,letterpaper]{article}
\usepackage{algorithm}
\usepackage{algorithmic}
\usepackage{subfigure}
\usepackage{amsmath}
\usepackage{mathtools}
\usepackage{multirow}
\usepackage{float}
\usepackage{booktabs}
\usepackage{graphicx}
\usepackage{amssymb}
\usepackage{stfloats}

\usepackage[top=0.85in,left=2.75in,footskip=0.75in,marginparwidth=2in]{geometry}

\usepackage[utf8]{inputenc}

\usepackage{cite}

\usepackage{nameref,hyperref}

\usepackage[right]{lineno}

\usepackage{microtype}
\DisableLigatures[f]{encoding = *, family = * }

\raggedright
\usepackage{changepage}

\usepackage[aboveskip=1pt,labelfont=bf,labelsep=period,singlelinecheck=off]{caption}

\makeatletter
\renewcommand{\@biblabel}[1]{\quad#1.}
\makeatother

\usepackage{lastpage,fancyhdr,graphicx}
\usepackage{epstopdf}
\fancyheadoffset[L]{2.25in}
\fancyfootoffset[L]{2.25in}

\usepackage{color}

\definecolor{Gray}{gray}{.25}

\usepackage{graphicx}

\usepackage{sidecap}

\usepackage{wrapfig}
\usepackage[pscoord]{eso-pic}
\usepackage[fulladjust]{marginnote}
\reversemarginpar

\begin{document}
\bibliographystyle{unsrt}
\vspace*{0.35in}

\begin{flushleft}
{\Large
\textbf\newline{FCNCP: A Coupled Nonnegative CANDECOMP/PARAFAC Decomposition Based on Federated Learning}
}
\newline
\\
Yukai Cai\textsuperscript{1,2},
Hang Liu\textsuperscript{1,2,*},
Xiulin Wang\textsuperscript{3,4,*},
Hongjin Li\textsuperscript{1,2},
Ziyi Wang\textsuperscript{1,2},
Chuanshuai Yang\textsuperscript{1,2},
Fengyu Cong\textsuperscript{1,2,5,6}
\\
\bigskip
\bf{1} School of Biomedical Engineering, Faculty of Medicine, Dalian University of Technology, Dalian, China
\\
\bf{2} Key Laboratory of Integrated Circuit and Biomedical Electronic System, Liaoning Province, Dalian University of Technology, Dalian,China
\\
\bf{3} Stem Cell Clinical Research Center, the First Affiliated Hospital of Dalian Medical University, Dalian, China
\\
\bf{4} Dalian Innovation Institute of Stem Cell and Precision Medicine, Dalian, China
\\
\bf{5} Faculty of Information Technology, University of Jyv$\ddot{a}$skyl$\ddot{a}$,Jyv$\ddot{a}$skyl$\ddot{a}$, Finland
\\
\bf{6} Key Laboratory of Social Computing and Cognitive Intelligence (Dalian University of Technology), Ministry of Education, China
\\
\bigskip
* liuhang@dlut.edu.cn
\\
* xiulin.wang@foxmail.com

\end{flushleft}

\section*{Abstract}
In the field of brain science, data sharing across servers is becoming increasingly challenging due to issues such as industry competition, privacy security, and administrative procedure policies and regulations. Therefore, there is an urgent need to develop new methods for data analysis and processing that enable scientific collaboration without data sharing. In view of this, this study proposes to study and develop a series of efficient non-negative coupled tensor decomposition algorithm frameworks based on federated learning called FCNCP for the EEG data arranged on different servers. It combining the good discriminative performance of tensor decomposition in high-dimensional data representation and decomposition, the advantages of coupled tensor decomposition in cross-sample tensor data analysis, and the features of federated learning for joint modelling in distributed servers. The algorithm utilises federation learning to establish coupling constraints for data distributed across different servers. In the experiments, firstly, simulation experiments are carried out using simulated data, and stable and consistent decomposition results are obtained, which verify the effectiveness of the proposed algorithms in this study. Then the FCNCP algorithm was utilised to decompose the fifth-order event-related potential (ERP) tensor data collected by applying proprioceptive stimuli on the left and right hands. It was found that contralateral stimulation induced more symmetrical components in the activation areas of the left and right hemispheres. The conclusions drawn are consistent with the interpretations of related studies in cognitive neuroscience, demonstrating that the method can efficiently process higher-order EEG data and that some key hidden information can be preserved.

\newcommand\keywords[1]{\textbf{Keywords}: #1}
 
\keywords{Non-negative coupled tensor decomposition, Federated learning, Event-related potentials (ERP)}

\section*{Introduction}
Hitchcock initially introduced the concept of tensors in the literature\cite{1927The}. In simple terms, tensors are multidimensional arrays, representing an extension of data representation from vectors and matrices to higher-dimensional spaces. They provide an intuitive means to characterize the structural features of high-dimensional data. From the standpoint of signal processing and data analysis, tensor decomposition allows for the consideration of information in multiple dimensions, encompassing time, space, and frequency. The resulting factor matrices or hidden components often possess physical and physiological significance. As a result, tensor decomposition has found extensive applications in various fields, including signal processing, machine learning, and has gained increasing prominence in the domains of brain signal processing and cognitive neuroscience\cite{2016Tensor,2015Tensor,2015Linked,2011Applications}. At the model level, two of the most widely used tensor decomposition models are the Canonical Polyadic (CP) model\cite{1927The}, which decomposes tensors into a sum of rank-one components, and the Tucker model\cite{1966Some}, which breaks down tensors into a kernel tensor and a series of factor matrices. Both of these models impose strict multicollinearity assumptions, which may not hold true in practical data processing scenarios. The PARAFAC2 model\cite{Henk1999PARAFAC2} addresses this issue by introducing an evolving factor, allowing tensor data to change over time within certain factor matrices. This enables the tracking of temporal changes in brain connectivity networks\cite{2020Tracing}. From a methodological perspective, the classical Alternating Least Squares (ALS) algorithm leverages an alternating least squares optimization criterion for parameter estimation\cite{2009Nonnegative}. The Hierarchical Alternating Least Squares (HALS) algorithm transforms the ALS optimization problem into multiple local optimization tasks, making it suitable for overdetermined or underdetermined cases. The Fast HALS algorithm enhances computational efficiency while maintaining separation accuracy through mathematical optimization of the HALS process\cite{cichocki2007hierarchical,cichocki2009fast}. Additionally, optimization algorithms based on the Effective Set Method, Alternating Proximal Gradient, and Alternating Direction Multiplier Methods have demonstrated exceptional effectiveness in tensor decomposition\cite{boyd2011distributed,kim2012fast,xu2015alternating,wang2021inexact,2014Lowrank}.

Coupled tensor decomposition is a method proposed in recent years for the joint analysis of multiple tensor data, extending the concept of tensor decomposition to address the joint blind source separation problem\cite{2016Linked}. Coupled tensor decomposition, when applied to tensor data from various samples, leverages its high-order and multiple coupling properties to facilitate joint analysis. This approach extracts a richer set of commonalities and differences between samples, thus enriching our understanding of brain function. It reveals potential internal mechanisms underpinning the relationship between brain function and various cognitive tasks or brain diseases. Numerous studies have explored joint analyses of EEG and fMRI data using coupled matrix-tensor factorization (CMTF\cite{acar2013structure}) and its variations\cite{2019Unraveling, 2020Group, 2022Exploring}. In these contexts, EEG data is typically represented as a third-order tensor, while fMRI data takes the form of a second-order matrix. These studies often assume shared structures across time or subject dimensions\cite{hunyadi2017tensor,martinez2004concurrent,acar2017tensor,jonmohamadi2020extraction}. For instance, Acar et al. conducted a joint analysis of EEG and fMRI data within an auditory oddball paradigm using the CMTF algorithm. In this study, EEG data is represented as a third-order tensor involving subjects, time, and channels, while fMRI data is presented as a second-order matrix with subjects and voxels. This approach successfully identified biologically significant components that could discriminate between patients with schizophrenia and healthy individuals. It achieved this through complementary high spatiotemporal resolution, thereby improving the performance of component clustering compared to using only the EEG tensor\cite{acar2017tensor}. However, it's important to note that existing EEG data processing and analysis techniques are primarily tailored for aggregated data on a single server. They lack the capacity to handle ``siloed'' data distributed across different servers and other scenarios.

In the era of advancing machine learning and big data technologies, various challenges have emerged, including concerns related to data authentication, privacy protection, and adherence to local policies and regulations. These challenges have led to the emergence of ``data silos'', where data circulation is hindered. Traditional approaches like edge computing and software-level data encryption are inadequate for effective coordination among multiple parties for collaborative training and learning\cite{shi2016edge,yu2020new}. To address these issues, researchers have proposed innovative machine learning frameworks based on federated learning\cite{li2020review,li2020federated,zhang2021survey,mammen2021federated,yang2019federated}. Unlike conventional centralized learning, which local servers send data to a central server for unified training, or distributed learning approaches like edge computing, which local servers train independently and then share results with a central server for aggregation and storage, federated learning adopts a local computation model transfer strategy. In this approach, servers do not exchange original data; instead, they transmit intermediate parameter results from local data training models to the central server for aggregation. Simultaneously, the aggregated parameters are sent back to local servers for individual updates, creating an iterative process that continues until the completion of training. This approach effectively involves local servers in joint model building while safeguarding the privacy and security of original data\cite{zhang2021survey,abdulrahman2020survey,mcmahan2017communication}. Depending on the data's feature space and sample space distribution, federated learning can be categorized into horizontal federated learning, vertical federated learning, and federated migration learning\cite{zhang2021survey}.

Federated tensor decomposition is a recently introduced model, combining tensor decomposition with federated learning. In this model, data residing on local servers possesses a tensor structure. Through iterative tensor decomposition, local factor matrices are computed, uploaded to the central server, and aggregated to generate the global factor matrix, which is then distributed. Each local server receives this global factor matrix as an initialization to commence the subsequent iteration updates, continuing until the iteration process concludes, resulting in the final global factor matrix\cite{kong2019federated,kim2017federated,gao2021federated}. Kim et al. achieved a significant milestone by using the federated tensor decomposition model to phenotype a vast Electronic Health Record (EHR) dataset spanning multiple hospitals. Crucially, this was accomplished without sharing patient-level data. The data from each hospital can be represented as a third-order tensor encompassing patients, medications, and diagnoses. Federated tensor decomposition proved its effectiveness by facilitating secure data coordination through the alternating direction method of multipliers and joint computation steps. Importantly, it matched or exceeded the accuracy and phenotype mining performance of centrally trained models while preserving patient privacy\cite{kim2017federated}. Additionally, some researchers proposed a simultaneous tensor decomposition model within the federated learning framework. They introduced a joint higher-order orthogonal iteration method to optimize the model, enabling secure joint data computation, feature extraction, and dimensionality reduction. Their work demonstrated that the federated tensor decomposition model can ensure both decomposition accuracy and classification accuracy while safeguarding data privacy, using real industrial datasets\cite{gao2021federated}. However, it's worth noting that current federated tensor decomposition models are designed exclusively for cross-server tensor decomposition. This implies an assumption that the latent component information in data from various servers is uniform. As of now, no research has explored the integration of coupled tensor decomposition models with the federated learning framework to extract information regarding the similarity and dissimilarity of data from distinct servers.

In our study, we proposed FCNCP, a coupled tensor decomposition model based on a federated framework, which not only retains the advantages of coupled tensor decomposition in cross-sample tensor data analysis, but also solves the problem that the existing EEG data processing and analysis techniques can only be used to process the aggregated data on a single server, and are not able to analyse ``siloed'' data from different servers and so on. This paper makes the following contributions.

1) We proposed FCNCP, a coupled tensor decomposition algorithm based on the federation framework. The algorithm uses federation learning to establish coupling constraints for coupled tensor decomposition. It is able to solve the problem of achieving joint data analysis without data sharing (privacy protection).

2) In the algorithmic flow of this model, a method for selecting coupling components (local public components) based on correlation analysis is proposed. This method can improve the stability of the model results and is a necessary part of the process of establishing coupling constraints through federated learning.

3) In order to demonstrate the effectiveness of the method, we designed simulation experiments using simulated data for third-order FCNCP decomposition, and concluded that the method is feasible and effective through the consistency and stability of the experimental results.

4) We apply the algorithm to ERP real data with proprioceptive stimuli applied to the right and left hands, and by comparing it with existing conclusions, we reach conclusions that are consistent with the interpretations of relevant studies in cognitive neuroscience, thus demonstrating that the method can effectively process higher-order EEG data and that some key hidden information can be preserved.

\section*{Method}

\subsection*{Overview of the proposed architecture} 

\begin{figure}[ht] 
\includegraphics[width=\textwidth]{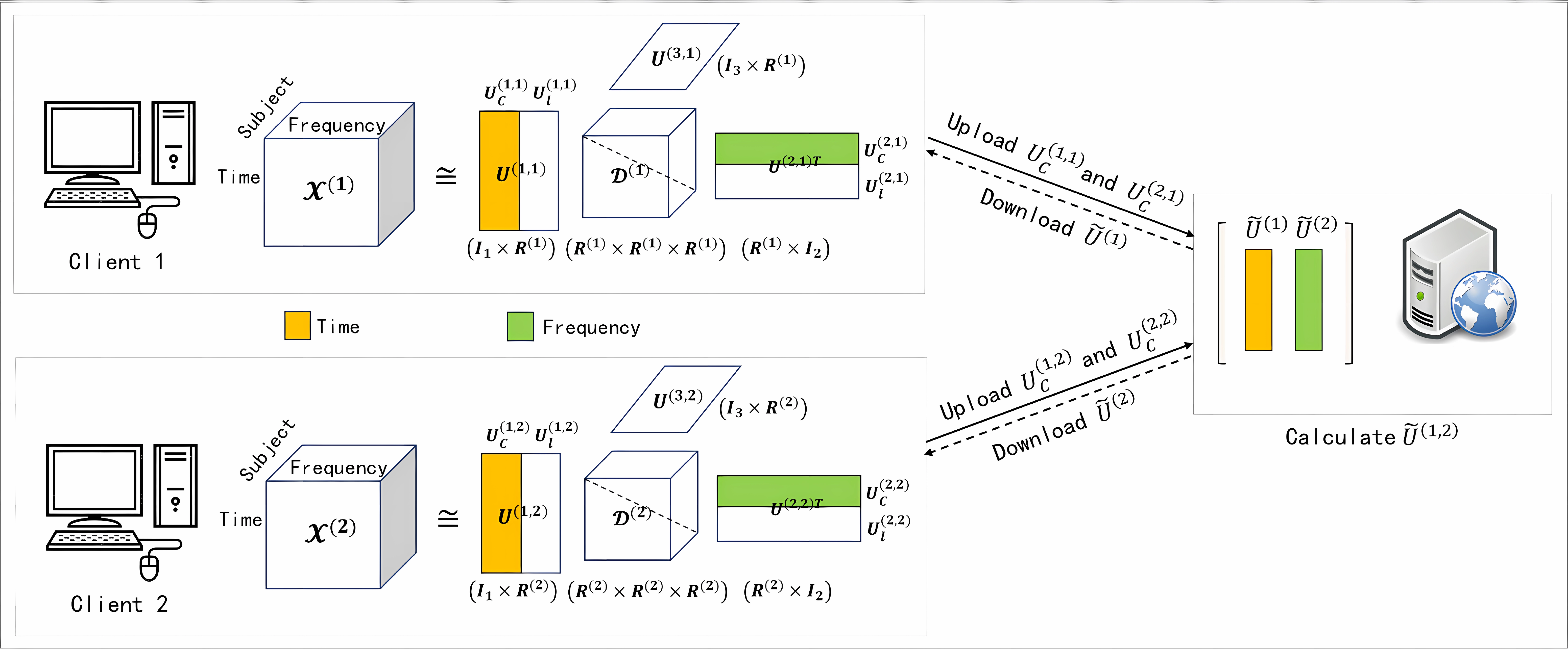}

\caption{Schematic diagram of the coupled tensor decomposition framework based on federated learning.}

\label{fig1} 

\end{figure}

In our study, we first represent EEG data within each client using a high-dimensional tensor format, considering the distribution characteristics of the data and its application scenarios. This tensor is structured as a third-order tensor, incorporating the dimensions of time, frequency, and subjects (trial times), as depicted in Figure~\ref{fig1}. Employing tensor decomposition, we break down the EEG data into three factor matrices that correspond to the dimensions of time, frequency, and subjects. Furthermore, we make an assumption regarding the sharing of temporal and frequency information among subjects within client groups. This assumption is visually represented by the yellow and green boxes in the figure. The remaining portion of the data comprises temporal and frequency information that is privately owned by subjects within client groups. Under the federated learning framework, the computation and updating of the latent factor matrix are divided into three distinct steps. This structured approach enables us to better understand the underlying factors in EEG data and its distribution.

Step 1: The tensor decomposition is run for the $i$th iteration for the data on clients 1 and 2, respectively, and the local public components representing the shared temporal and frequency information in the iteration results are uploaded to the central server, respectively. It is important to note that the constituents for which coupling constraints are to be established need to be selected in the central server at the 1st iteration and their positions returned to the client.

Step 2: The central server calculates global model parameter by fusing the local public components from these clients and sends them down to clients 1 and 2.

Step 3: Clients 1 and 2 update the local public components after receiving the global model parameter and use them as initial variables for the $i+1$ th iteration.

Repeat the above three steps until the local public components satisfy the iteration stop condition.


\subsection*{Basic symbols and mathematical operations}
${{\left\| \centerdot  \right\|}_{F}}$ and ${{\left( \centerdot  \right)}^{T}}$ denote  the Frobenius paradigm of the matrix and the matrix transpose, respectively, and the outer product, Khatri-Rao product and Hadamard product are denoted by '$\circ$', '$\odot$', and '$\circledast$'. In addition, ${{U}^{\left( 1 \right)}}\odot {{U}^{\left( 2 \right)}}\odot \cdots \odot {{U}^{\left( N \right)}}$, ${{U}^{\left( 1 \right)}}\odot \cdots \odot {{U}^{\left( n-1 \right)}}\odot {{U}^{\left( n+1 \right)}}\odot \cdots \odot {{U}^{\left( N \right)}}$, ${{U}^{\left( 1 \right)}}\circledast {{U}^{\left( 2 \right)}}\circledast \cdots \circledast {{U}^{\left( N \right)}}$, and ${{U}^{\left( 1 \right)}}\circledast \cdots \circledast {{U}^{\left( n-1 \right)}}\circledast {{U}^{\left( n+1 \right)}}\circledast \cdots \circledast {{U}^{\left( N \right)}}$ are defined by ${{\left\{ U \right\}}^{\odot }}$, ${{\left\{ U \right\}}^{{{\odot }_{-n}}}}$, ${{\left\{ U \right\}}^{\circledast }}$, and ${{\left\{ U \right\}}^{{{\circledast }_{-n}}}}$, respectively. The mode-n expansion of a tensor $\chi \in {{\mathbb{R}}^{{{I}_{1}}\times {{I}_{2}}\times \cdots \times {{I}_{N}}}}$ is defined as ${{\chi }_{\left( n \right)}}$, and its size is ${{I}_{n}}\times \left( {{I}_{1}}\cdots {{I}_{n-1}}{{I}_{n+1}}\cdots {{I}_{N}} \right)$.


\subsection*{Coupled tensor decomposition}
Coupled tensor decomposition is a method for simultaneous decomposition of multiple tensors. While traditional tensor decomposition methods can only handle a single tensor, coupled tensor decomposition can efficiently handle multiple related tensors. In coupled tensor decomposition, we consider the correlations between multiple tensors and use these correlations to improve the accuracy and efficiency of the decomposition.

Assuming that there are $S$ tensor blocks ${{\chi }^{\left( s \right)}}$, the corresponding CP decomposition is denoted as ${{\chi }^{\left( s \right)}}\simeq \sum\nolimits_{r=1}^{{{R}^{\left( s \right)}}}{u_{r}^{\left( 1,s \right)}\circ u_{r}^{\left( 2,s \right)}\circ \cdots \circ u_{r}^{\left( N,s \right)}}$, and the objective function is constructed considering the non-negativity of the constructed tensor and the coupling constraints in the factor matrix:

\begin{align}
    \begin{split}
        & \min \sum\limits_{s=1}^{S}{\left\| {{\chi }^{\left( s \right)}}-\sum\nolimits_{r=1}^{{{R}^{\left( s \right)}}}{u_{r}^{\left( 1,s \right)}\circ u_{r}^{\left( 2,s \right)}\circ \cdots \circ u_{r}^{\left( N,s \right)}} \right\|_{F}^{2}} \\ 
        & s.t.\text{ }u_{r}^{\left( n,1 \right)}=u_{r}^{\left( n,2 \right)}=\cdots =u_{r}^{\left( n,S \right)},r\le {{L}_{n}} \\ 
    \end{split}
    \label{equation1}
\end{align}

where $N$ denotes the number of tensor dimensions, ${{R}^{\left( s \right)}}$ denotes the number of components of the $s$th tensor block, and $S$ denotes the number of tensor blocks.

\subsection*{Federated learning}
The goal of federated learning algorithms is to learn a robust model from multicentric data. Federated Averaging (FedAvg) is a simultaneous distributed optimisation algorithm which is the best known FL algorithm\cite{mcmahan2017communication}. It consists of local updating and global aggregation. Its learning process can be described by (\ref{equation2}).

\begin{align}
    \begin{split}
        &w_{t}^{k}\triangleq
    {\begin{cases}
             w_{t-1}^{k}-\eta\nabla{{F}^{k}}\left( w_{t-1}^{k} \right),&if \text{ }t \text{ } mod \text{ }\tau\ne0\\
             {{w}_{t}},&if \text{ }t \text{ } mod \text{ }\tau=0
    \end{cases}}\\
    &where\text{ }{{w}_{t}}\triangleq \sum\limits_{k\in N}{\frac{\left| {{D}^{k}} \right|}{\left| D \right|}\left[ w_{t-1}^{k}-\eta \nabla {{F}^{k}}\left( w_{t-1}^{k} \right) \right]}
    \end{split}
    \label{equation2}
\end{align}

Where, $\left( x,y \right)\in {{D}^{k}}$ denotes the data set of the $k$th client, and $w_{t}^{k}$ denotes the parameter of the kth client in the tth round. The formula represents the process of gradient descent with learning rate $\eta$ and aggregating the parameters of the clients every $\tau$ rounds.

\subsection*{A coupled tensor decomposition model based on an elastic average federal framework}
In model averaging algorithms for federated learning, a global model is aggregated at some point and used in place of the local model. However, this approach is not the most suitable choice for optimisation problems with many local minima, as the locally trained models will differ due to the different data of each client, and each client actually searches for the local optimal point in a different search space, and due to the different directions of exploration, it is possible to get very different models.

In order to solve the above problem, researchers proposed a non-exactly consistent federated learning algorithm called Elastic Mean Federated Learning Algorithm\cite{zhang2015deep}. The method does not force individual clients to inherit the global model (also the idea of personalised federated learning in later federated learning). We apply its idea to the coupled tensor decomposition model by using federated learning to establish coupling constraints, thus proposing a new federated learning-based coupled tensor decomposition model called FCNCP and constructing the following objective function:

\begin{align}
    \begin{split}
    & \min L=\sum\limits_{k=1}^{K}{\left\| {{\chi }_{k}}-\sum\nolimits_{r=1}^{{{R}_{k}}}{u_{r,k}^{\left( 1 \right)}\circ u_{r,k}^{\left( 2 \right)}\circ \cdots \circ u_{r.k}^{\left( N \right)}} \right\|_{F}^{2}}\\
    &+\frac{\rho }{2}\sum\limits_{k=1}^{K}{\sum\limits_{n=1}^{N}{\sum\limits_{r=1}^{{{L}_{n}}}{\left\| u_{r,k}^{\left( n \right)}-\tilde{u}_{r}^{\left( n \right)} \right\|_{2}^{2}}}} \\ 
    & s.t.\text{ }u_{r,k}^{\left( n \right)}\ge 0\text{ }for\text{ }n=1,2,\cdots N
    \end{split}
    \label{equation3}
\end{align}

where $\tilde{u}_{r}^{\left( n \right)}$ is the global model and $K$ is the number of clients. This distributed optimisation method has two objectives, one is to make the loss function local to each client to be minimised, and the second is that the gap between the local model on each client and the global model on the central server is desired to be relatively small.

\subsection*{Federal optimisation methodology}
Our goal is to find local model parameters $u_{r,k}^{\left( n \right)}$ and global model parameters $\tilde{u}_{r}^{\left( n \right)}$ to minimise the objective (\ref{equation3}). In this paper, an elasticity-averaging based fast hierarchical alternating least squares algorithm (EAFHALS) is used.

\subsubsection*{Local personalised model update}

Each client updates the local model (including local private and public components) by solving the subproblems of (\ref{equation3}).

Local private components update: Different clients have different local private components, i.e., $u_{r,k}^{\left( n \right)},r>{{L}_{n}}$, and using the derivation of (\ref{equation3}) to obtain the gradient equation as follows:

\begin{equation}
    \frac{\partial L}{\partial u_{r,k}^{\left( n \right)}}=-2{\mathcal{Y}_{r,\left( n \right)}}{{\left\{ {{u}_{r,k}} \right\}}^{{{\odot }_{-n}}}}+2u_{r,k}^{\left( n \right)}{{\left\{ u_{r,k}^{T}{{u}_{r,k}} \right\}}^{{{\circledast }_{-n}}}}\
\end{equation}

where ${\mathcal{Y}_{r,k}}={{\chi }_{k}}-\sum\nolimits_{k\ne r}^{{{R}_{k}}}{u_{r,k}^{\left( 1 \right)}\circ u_{r,k}^{\left( 2 \right)}\circ \cdots \circ u_{r.k}^{\left( N \right)}}$. For each client, the gradient equation is made equal to 0. The client's local private component is updated according to its local dataset, and the update equation is as follows:

\begin{equation}
    u_{r,k}^{\left( n \right),t+1}=u_{r,k}^{\left( n \right),t}+\frac{{{\left[ {{\chi }_{k\left( n \right)}}{{\left\{ {{U}_{k}} \right\}}^{{{\odot }_{-n}}}} \right]}_{r}}-U_{k}^{\left( n \right)}{{\left[ {{\left\{ U_{k}^{T}{{U}_{k}} \right\}}^{{{\circledast }_{-n}}}} \right]}_{r}}}{{{\left\{ u_{r,k}^{T}{{u}_{r,k}} \right\}}^{{{\circledast }_{-n}}}}}\
    \label{equation5}
\end{equation}

Local public component update: The local public components are the coupled component in the coupled tensor decomposition, i.e., $u_{r,k}^{\left( n \right)},r\le {{L}_{n}}$, and using the derivation of (\ref{equation3}) to obtain the gradient equation as follows:

\begin{align}
    \begin{split}
        \frac{\partial L}{\partial u_{r,k}^{\left( n \right)}}=&-2{\mathcal{Y}_{r,\left( n \right)}}{{\left\{ {{u}_{r,k}} \right\}}^{{{\odot }_{-n}}}}+2u_{r,k}^{\left( n \right)}{{\left\{ u_{r,k}^{T}{{u}_{r,k}} \right\}}^{{{\circledast }_{-n}}}}\\
        &+\rho \left( u_{r,k}^{\left( n \right)}-\tilde{u}_{r}^{\left( n \right)} \right)\
    \end{split}
\end{align}

For each client, the gradient equation is made equal to 0. The client's local public component is updated according to its local dataset, and the update equation is as follows:

\begin{align}
    \begin{split}
        u_{r,k}^{\left( n \right),t+1}=&u_{r,k}^{\left( n \right),t}\cdot \frac{{{\left\{ u_{r,k}^{T}{{u}_{r,k}} \right\}}^{{{\circledast }_{-n}}}}}{{{\left\{ u_{r,k}^{T}{{u}_{r,k}} \right\}}^{{{\circledast }_{-n}}}}+\frac{\rho }{2}}+\\
        &\frac{{{\left[ {{\chi }_{k\left( n \right)}}{{\left\{ {{U}_{k}} \right\}}^{{{\odot }_{-n}}}} \right]}_{r}}-U_{k}^{\left( n \right)}{{\left[ {{\left\{ U_{k}^{T}{{U}_{k}} \right\}}^{{{\circledast }_{-n}}}} \right]}_{r}}+\frac{\rho }{2}\tilde{u}_{r}^{\left( n \right)}}{{{\left\{ u_{r,k}^{T}{{u}_{r,k}} \right\}}^{{{\circledast }_{-n}}}}+\frac{\rho }{2}}
    \end{split}
    \label{equation7}
\end{align}

\subsubsection*{Global model parameter update}

The central server accepts the local public components from each client $u_{r,k}^{\left( n \right)},r\le {{L}_{n}}$, and calculates the global model parameters $\tilde{u}_{r}^{\left( n \right)}$ based on the local public components, and broadcasts it to each client, and the gradient equation can be obtained by taking the partial derivation of the global model parameters from the (\ref{equation3}) as follows:

\begin{equation}
    \frac{\partial L}{\partial \tilde{u}_{r}^{\left( n \right)}}=\rho \left( \tilde{u}_{r}^{\left( n \right)}-\frac{1}{K}\sum\limits_{k=1}^{K}{u_{r,k}^{\left( n \right)}} \right),r\le {{L}_{n}}\
    \label{equation8}
\end{equation}

Here we update the global model parameters by the central server according to the update equation of stochastic gradient descent method as:

\begin{equation}
    \tilde{u}_{r}^{\left( n \right),t+1}=\tilde{u}_{r}^{\left( n \right),t}-\alpha \frac{\partial L}{\partial \tilde{u}_{r}^{\left( n \right)}},r\le {{L}_{n}}\
    \label{equation9}
\end{equation}\\
where $\alpha$ is the learning rate.

In order to speed up the algorithm when iterating, we used the matrix of factors generated by the tensor block ${{\chi }_{k}}$ from its unconstrained CP decomposition $\left[\! \left[ \widehat{U}_{k}^{\left( 1 \right)},\widehat{U}_{k}^{\left( 2 \right)},\cdots ,\widehat{U}_{k}^{\left( N \right)} \right]\! \right]$ for low-rank approximation. Thus, the n-mode expansion is given by ${{\chi }_{k(n)}}=\widehat{U}_{k}^{\left( n \right)}{{\left[ {{\left\{ \widehat{U}_{k}^{\left( 1 \right)} \right\}}^{{{\odot }_ {-n}}}} \right]}^{T}}$ denotes.

\subsubsection*{Coupling component selection}
In the initial iteration of this federated algorithm, the selection of components requiring coupling is a critical step. Following the first round of decomposition for both clients' data, the factor matrices that necessitate the establishment of coupling constraints are uploaded. These matrices undergo computation of correlation coefficient matrices on the central server, followed by summation. The central server needs to identify and record the row and column indices associated with the ${{L}_{n}}$ highest values within these matrices. These indices are then transmitted back to clients 1 and 2. Subsequently, clients 1 and 2 categorize the components that correspond to the rows and columns where these indices are located as locally public components. At the same time, they designate the remaining components as locally private components. The specific algorithmic flow is shown in Algorithm \ref{algorithm2}.

\renewcommand{\algorithmicrequire}{ \textbf{Input:}} 
\renewcommand{\algorithmicensure}{ \textbf{Output:}} 

\begin{algorithm*}[htb]
    \caption{Coupled tensor decomposition federated learning optimisation algorithm}
    \label{algorithm1}
    \begin{algorithmic}[1]
        \REQUIRE ~~\\
        $u_{r,k}^{\left( n \right),t}\left( t=0,n=\left\{ 1,2,\cdots ,N \right\},r=\left\{ 1,2,\cdots ,{{R}_{k}} \right\},k=\left\{ 1,2, \right\} \right)$\\
        $\tilde{u}_{r}^{\left( n \right),t}\left( t=0,n=\left\{ 1,2,\cdots ,N \right\},r=\left\{ 1,2,\cdots ,{{L}_{n}} \right\} \right)$
        \ENSURE ~~\\
        $u_{r,k}^{\left( n \right)}\left( n=\left\{ 1,2,\cdots ,N \right\},r=\left\{ 1,2,\cdots ,{{R}_{k}} \right\},k=\left\{ 1,2 \right\} \right)$\\
        $\tilde{u}_{r}^{\left( n \right)}\left( n=\left\{ 1,2,\cdots ,N \right\},r=\left\{ 1,2,\cdots ,{{L}_{n}} \right\} \right)$
        
        \FOR{$t=0,1,\cdots,T-1$}
        \STATE
        Central server do
        \IF{$t=0$}
        \STATE
        Selection of coupling components using algorithm \ref{algorithm2}
        \ENDIF
        \STATE
        Receive local public components sent by each client\\
        $u_{r,k}^{\left( n \right)}\left( n=\left\{ 1,2,\cdots ,N \right\},r=\left\{ 1,2,\cdots ,{{L}_{n}} \right\},k=\left\{ 1,2 \right\} \right)$
        \STATE
        Normalising the uploaded coupling components
        \STATE
        Update the global model parameters $\tilde{u}_{r}^{\left( n \right)}\left( n=\left\{ 1,2,\cdots ,N \right\},r=\left\{ 1,2,\cdots ,{{L}_{n}} \right\} \right)$ according to (\ref{equation8}), (\ref{equation9})
        \STATE
        Broadcast global model parameters to individual clients
        \STATE
        
        \FOR{$k=1,2$}
        \STATE
        Client k do
        \IF{$t=0$}
        \STATE
        Labelling of local public components using algorithm \ref{algorithm2}
        \ENDIF
        \STATE
        Receive global model parameters sent by the central server $\tilde{u}_{r}^{\left( n \right)}\left( n=\left\{ 1,2,\cdots ,N \right\},r=\left\{ 1,2,\cdots ,{{L}_{n}} \right\} \right)$
        \STATE
        Update local private components $u_{r,k}^{\left( n \right)}\left( n=\left\{ 1,2,\cdots ,N \right\},r=\left\{ {{L}_{n}}+1,\cdots ,{{R}_{k}} \right\} \right)$ according to (\ref{equation5})
        \STATE
        Update the local public component $u_{r,k}^{\left( n \right)}\left( n=\left\{ 1,2,\cdots ,N \right\},r=\left\{ 1,2,\cdots ,{{L}_{n}} \right\} \right)$ according to (\ref{equation7})
        \STATE
        Upload local public composition to central server
        \ENDFOR
        \ENDFOR
    \end{algorithmic}
\end{algorithm*}

\renewcommand{\algorithmicrequire}{ \textbf{Input:}} 
\renewcommand{\algorithmicensure}{ \textbf{Output:}} 
\begin{algorithm}[htb]
\caption{Coupling component selection algorithm}
\label{algorithm2}
\begin{algorithmic}[1]
\REQUIRE ~~\\
$U_{k}^{\left( n \right )}\left(n=\left\{ 1,2,\cdots ,m \right\},m<N, k=\left\{1,2\right\}\right)$\\
$L_{n}\left(n=\left\{ 1,2,\cdots ,N \right\}\right)$
\ENSURE ~~\\
Component set $C_{k}\left(k=\left\{1,2\right\}\right)$
\STATE
Central server do
\FOR{$n=0,1,\cdots,m$}
\STATE
Calculate the matrix of correlation coefficients $P_{n}$ between $U_{n}^{1}$ and $U_{n}^{2}$
\ENDFOR
\STATE
$P=P_{1}+P_{2}+\cdots +P_{m}$
\STATE
$location_{1},location_{2} \gets emptyList$
\FOR{$l=1,2,\cdots max\left( L_{n}\right)$}
\STATE
Select the maximum value of the element in $P$ and record the row and column indexes $r$ and $c$ where it is located
\STATE
append($location_{1}$, $r$), append($location_{2}$, c)
\STATE
Place all elements of row $r$, column $c$ in 0
\ENDFOR
\STATE
Send $location_{1}$ and $location_{2}$ to client1 and 2 respectively
\STATE

\FOR{$k=1,2$}
\STATE
Client k do
\FOR{$n=1,2,\cdots N$}
\FOR{$l=1,2, \cdots, L_{n}$}
\STATE
Label the $location_{k}\left[l\right]$th component as the local public component of the $n$th pattern
\ENDFOR
\ENDFOR
\ENDFOR
\end{algorithmic}
\end{algorithm}

Algorithm \ref{algorithm1} describes the specific process of coupled tensor decomposition based on the elastic mean federal framework.

\section*{Experiments and results}
In this section, we first conducted simulation experiments using synthetic tensor data to verify the effectiveness of federated learning in establishing coupling constraints. Then, tensor decomposition experiments were conducted using the fifth-order tensor formed from the ERP data collected by applying proprioceptive stimuli to the right and left hands. By comparing the findings with the literature using the same dataset, it is demonstrated that the algorithm is able to process high-order EEG data efficiently with privacy preservation, and some key hidden information can be preserved.

\subsection*{Simulate experiment}
We first applied the algorithm with simulated data. Two tensors of size 61 × 72 × 64, representing frequency × time × channel, were created, and the two tensor blocks were used as two clients as follows:

\begin{align}
    \begin{split}
    &{{\chi }_{\text{1}}}={{t}_{1}}\circ {{f}_{1}}\circ {{c}_{1}}+{{t}_{2}}\circ {{f}_{2}}\circ {{c}_{2}}+{{t}_{3}}\circ {{f}_{3}}\circ {{c}_{3}}\\
    &{{\chi }_{2}}={{t}_{1}}\circ {{f}_{1}}\circ {{c}_{4}}+{{t}_{2}}\circ {{f}_{2}}\circ {{c}_{5}}+{{t}_{4}}\circ {{f}_{4}}\circ {{c}_{6}}
    \end{split}
\end{align}

\begin{figure*}[htb]
\centering

\vspace{-0.35cm}
\setlength{\abovecaptionskip}{-2pt}
\subfigtopskip=2pt
\subfigbottomskip=2pt
\subfigcapskip=-5pt

\subfigure[]{
\label{fig2a}
\includegraphics[width=6cm]{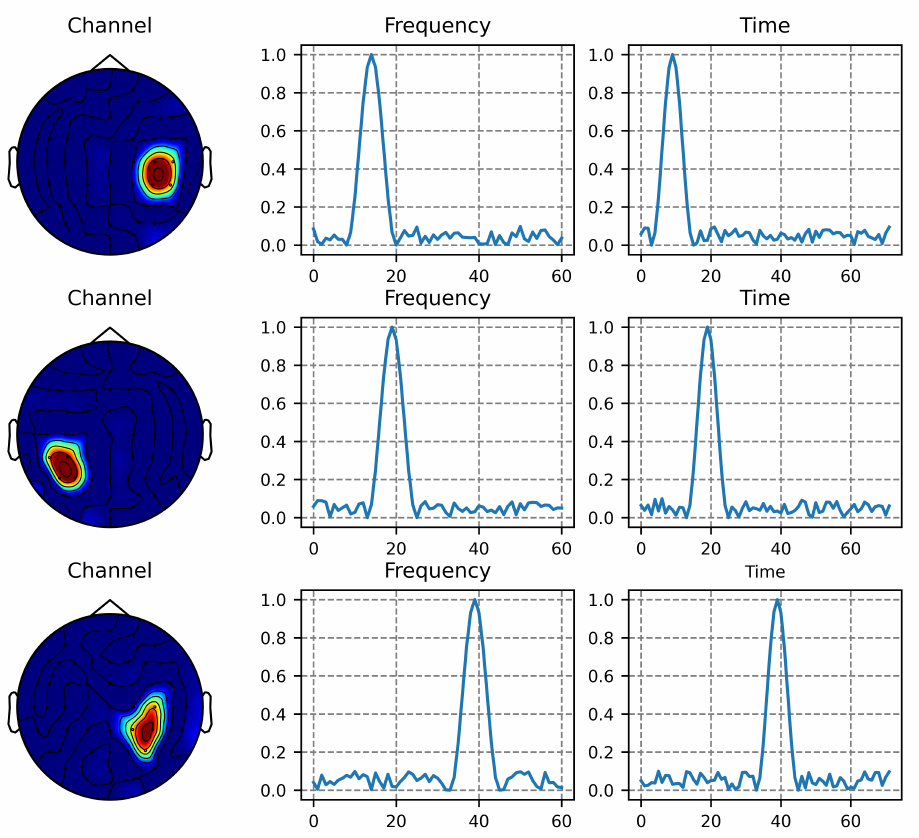}}
\quad
\subfigure[]{
\label{fig2b}
\includegraphics[width=6cm]{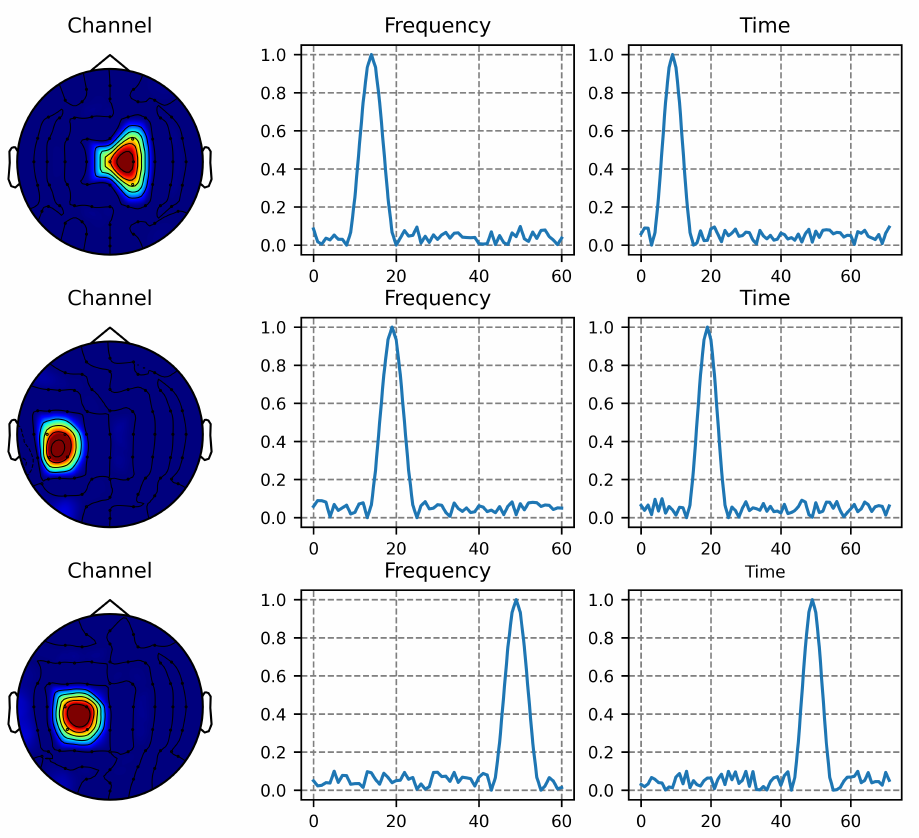}}
\subfigure[]{
\label{fig2c}
\includegraphics[width=6cm]{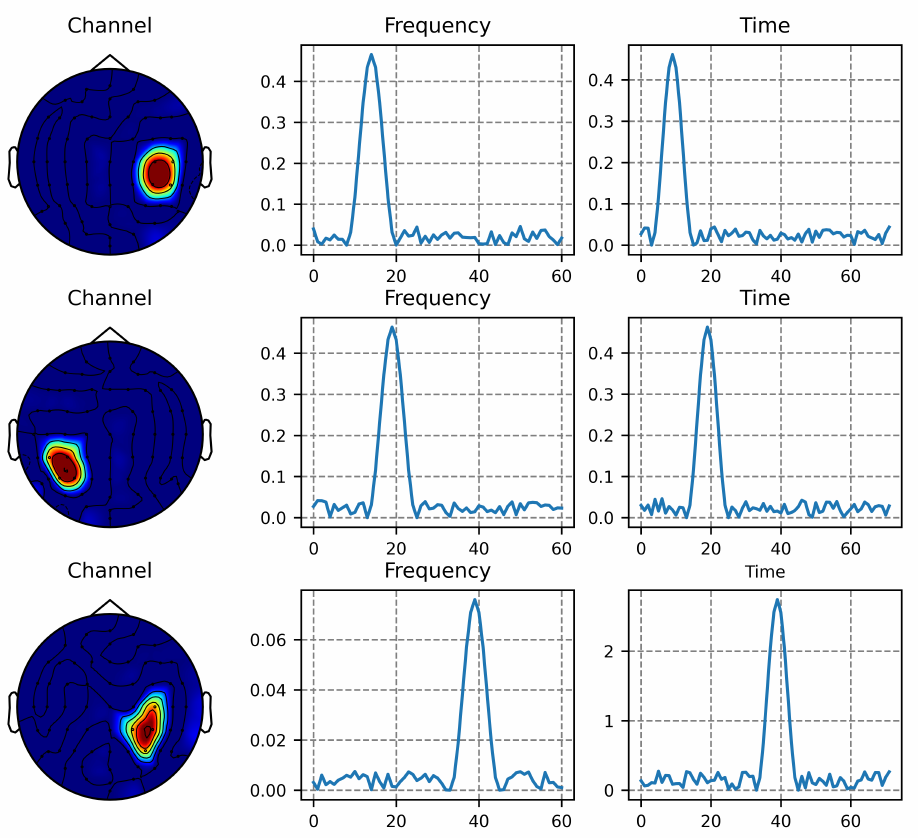}}
\quad
\subfigure[]{
\label{fig2d}
\includegraphics[width=6cm]{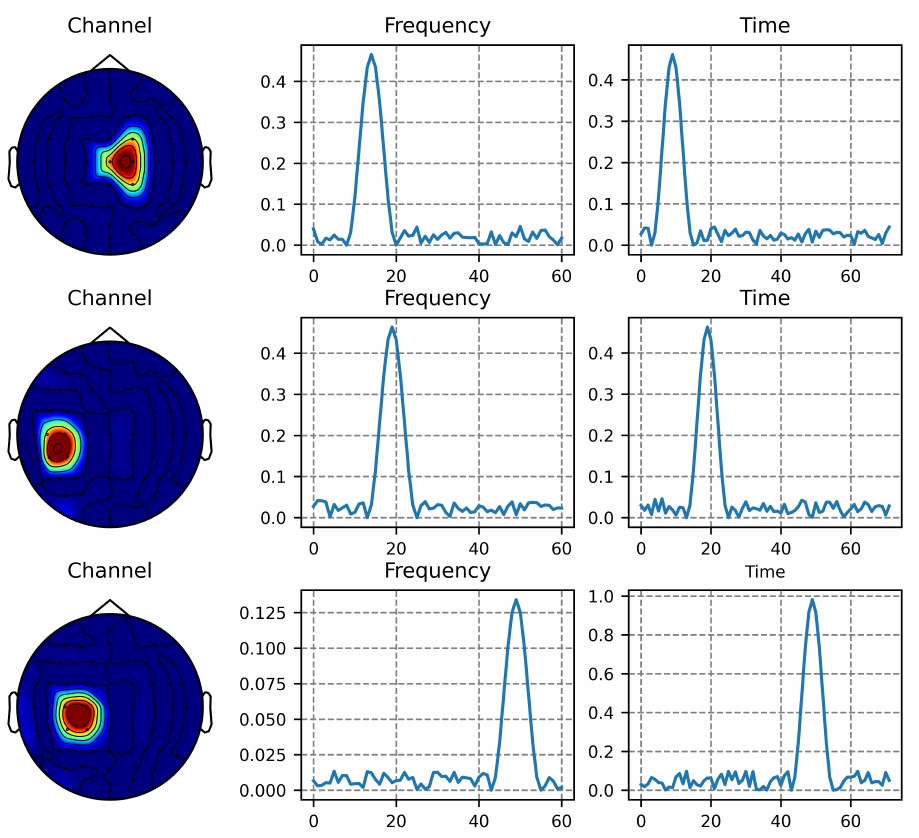}}
\caption{Illustration of the simulation experiment. (a) Tensor data for client 1. (b) Tensor data for client 2. (c) Results after FCNCP decomposition for client1. (d) Result after FCNCP decomposition for client2.}
\label{level}
\end{figure*}

In the frequency and time patterns, four frequency components were constructed using a Hanning window centred at 15, 20, 40, 50 and white noise, and four time components were constructed using a Hanning window centred at 10, 20, 30, 40 and white noise, and coupling was established on top of that. Six channel components were formed by randomly selecting four adjacent brain electrodes and assigning them a random number between 0.5 and 1. The simulated data for the two tensor blocks (client side) are shown in Figure~\ref{fig2a} and~\ref{fig2b}.

We ran 50 times of the FCNCP algorithm, and we obtained stable decomposition results with an averaged tensor fit of 0.996 for both client. According to the experimental results shown in Figure~\ref{fig2c} and~\ref{fig2d}, the coupling constraints we established in the simulation data can also be reflected at the final results after going through the decomposition of the FCNCP algorithm. This also verifies the effectiveness of federated learning in establishing coupling constraints.

\subsection*{Data description}
The ERP data in our experiments were obtained from the open preprocessing dataset associated with the ERPWAVELAB toolkit\cite{morup2007erpwavelab}, which can be downloaded from \href{www.erpwavelab.org}{www.erpwavelab.org}. These data were from proprioceptive experiments in which two conditions (left- and right-handed) were manipulated by increasing the handheld load. An important part of the stimulus was the change in force during static muscle contraction, which is considered a proprioceptive stimulus\cite{morup2006decomposing}. 14 subjects participated in the experiment and 64 scalp electrodes were used to record EEG data. Each subject performed a total of 360 trials (epochs) in each condition. All ephemeral data were transformed into time-frequency representation (TFR) by means of the complex Molette wavelet. In the wavelet transform, only the frequency band from 15 Hz to 75 Hz was analysed, with linear intervals of 1 Hz. Inter-trial phase coherence (ITPC) was then computed as an average spectral estimate over all trials\cite{delorme2004eeglab}. Since the TFR is first applied to each trial and then averaged across trials is calculated, the ITPC can be regarded as an induced oscillation of the brain\cite{david2006mechanisms}. Also, the ITPC only takes values between 0 and 1\cite{delorme2004eeglab,cohen2014analyzing}. Finally, a fifth-order non-negative tensor block (condition × subject × channel × frequency × time = 2 × 14 × 64 × 61 × 72) was generated, with 61 frequency points representing 15 to 75 Hz, and 72 time points representing 0 to 346.68 ms. Finally, a simulation experiment of federated learning was carried out in which the 14 subjects were divided into two clients, and assigned to seven subjects, respectively.

\subsection*{Experiment setup}\label{sec32}

\begin{figure}[htb]
    \centering
    \includegraphics[width=9cm]{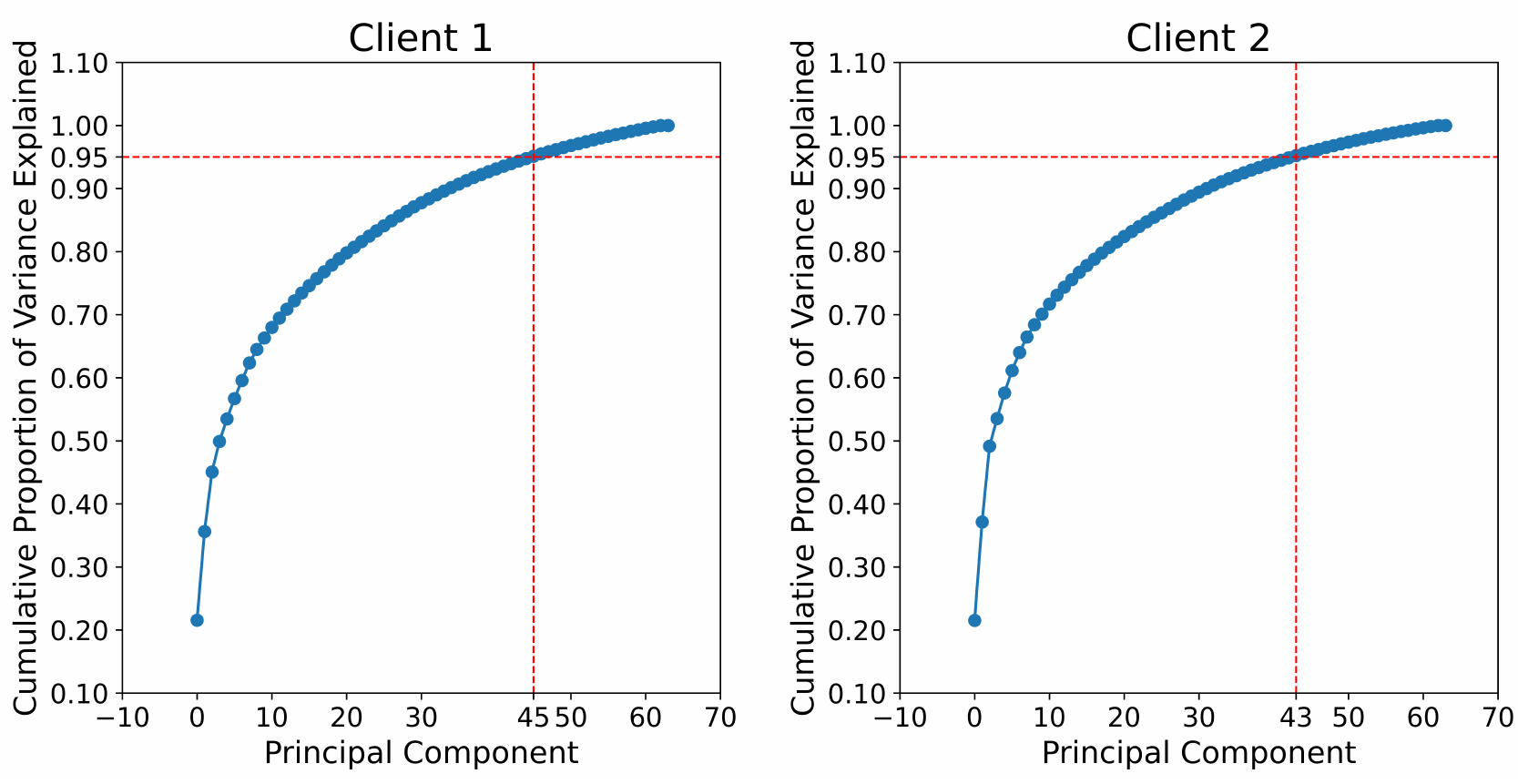}
    \caption{PCA was used to select the number of components with 95 percent explained cumulative variance.}
    \label{fig2}
\end{figure}

\begin{figure}[htb]
    \centering
    \includegraphics[width=11cm]{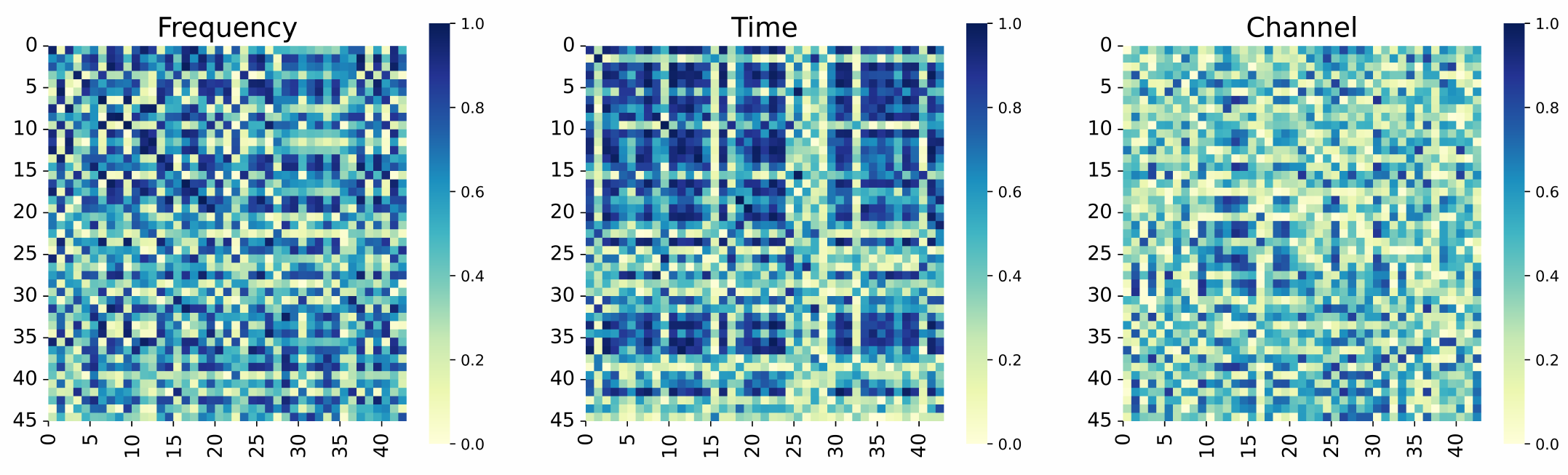}
    \caption{Correlation analysis of frequency factor, time factor and channel factor extracted by four-dimensional tensor decomposition for each block of data.}
    \label{fig3}
\end{figure}

In this section, we describe how to determine the number of components of each client tensor block when performing the coupled tensor decomposition ${{R}_{k}}\left( k=1,2 \right)$, i.e., the hidden information of each block of data in the low-dimensional space, and the number of coupled components. They reveal the common features between the blocks of data. For the selection of ${{R}_{k}}\left( k=1,2 \right)$ we perform a principal component analysis (PCA) on the matrix data of each chunk of data unfolded along the frequency pattern, and retain the number of components with 95\% cumulative variance explained, as in Figure~\ref{fig2}, where the two clients select \{45, 43\} as their number of components, respectively. The selection of the number of coupled components is a key issue for the coupled tensor decomposition algorithm to be performed and the results to be interpreted, and it is always an open problem depending on the practical application. In this paper, the fifth-order CP tensor decomposition based on the ALS algorithm is performed on the data in three clients respectively, and the factor matrices under the time, frequency and channel dimensions, which are needed to establish coupling constraints between the clients, are uploaded to the central server, and the correlation maps between the components extracted from the data in the time mode, the frequency mode, and the channel mode are computed respectively. Based on the correlation maps, we will select the number of highly correlated (coupled) components. As shown in Figure~\ref{fig3}, the tensor blocks in the three clients have high correlation in the time and frequency factors after fifth-order CP decomposition, while their correlation is low in the channel factor analysis, so we choose to establish the coupling constraints on time and frequency, and select 15 as the number of coupled components.

\begin{figure*}
    \centering
    \subfigtopskip=2pt
    \subfigbottomskip=2pt
    \subfigcapskip=-5pt
    \subfigure[]{\includegraphics[width=13cm]{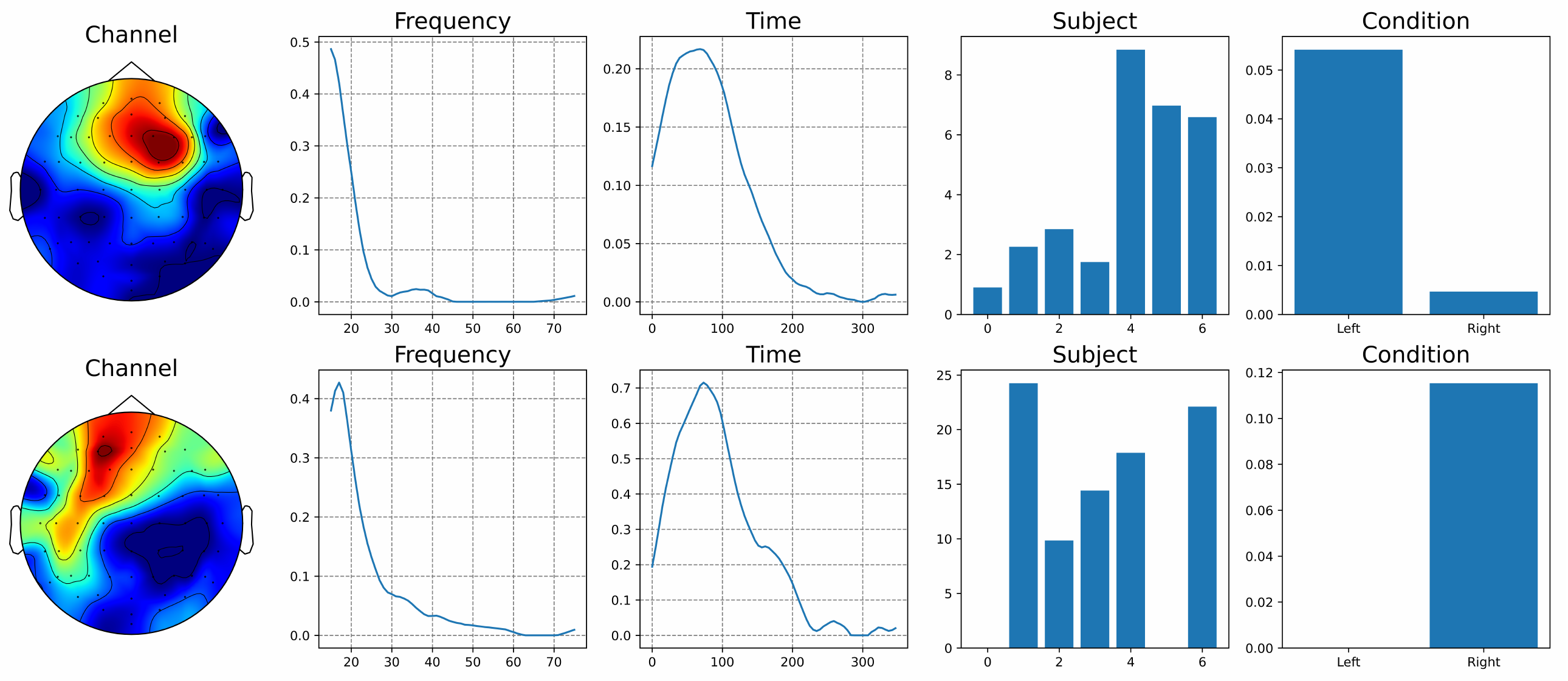}
    \label{fig5a}}
    \subfigure[]{\includegraphics[width=13cm]{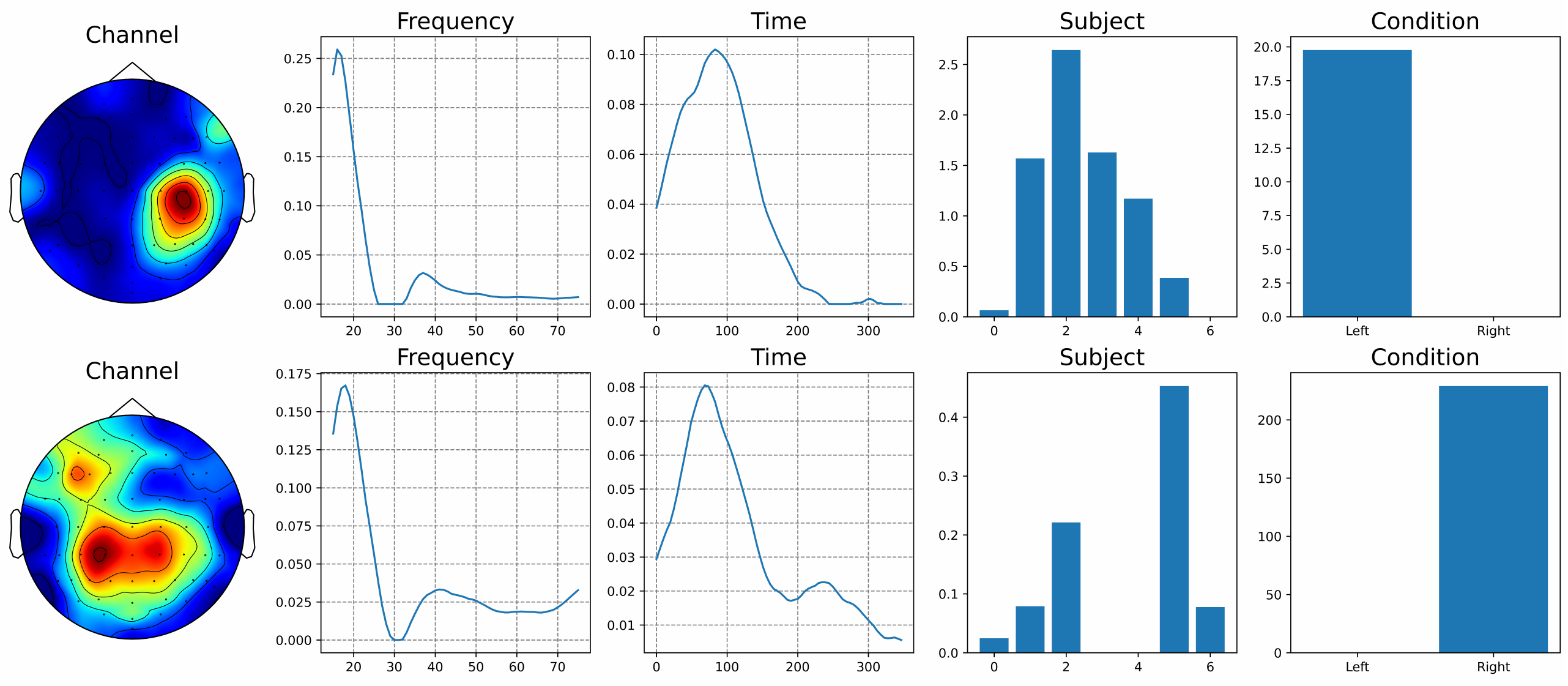}
    \label{fig5b}}
    \subfigure[]{\includegraphics[width=13cm]{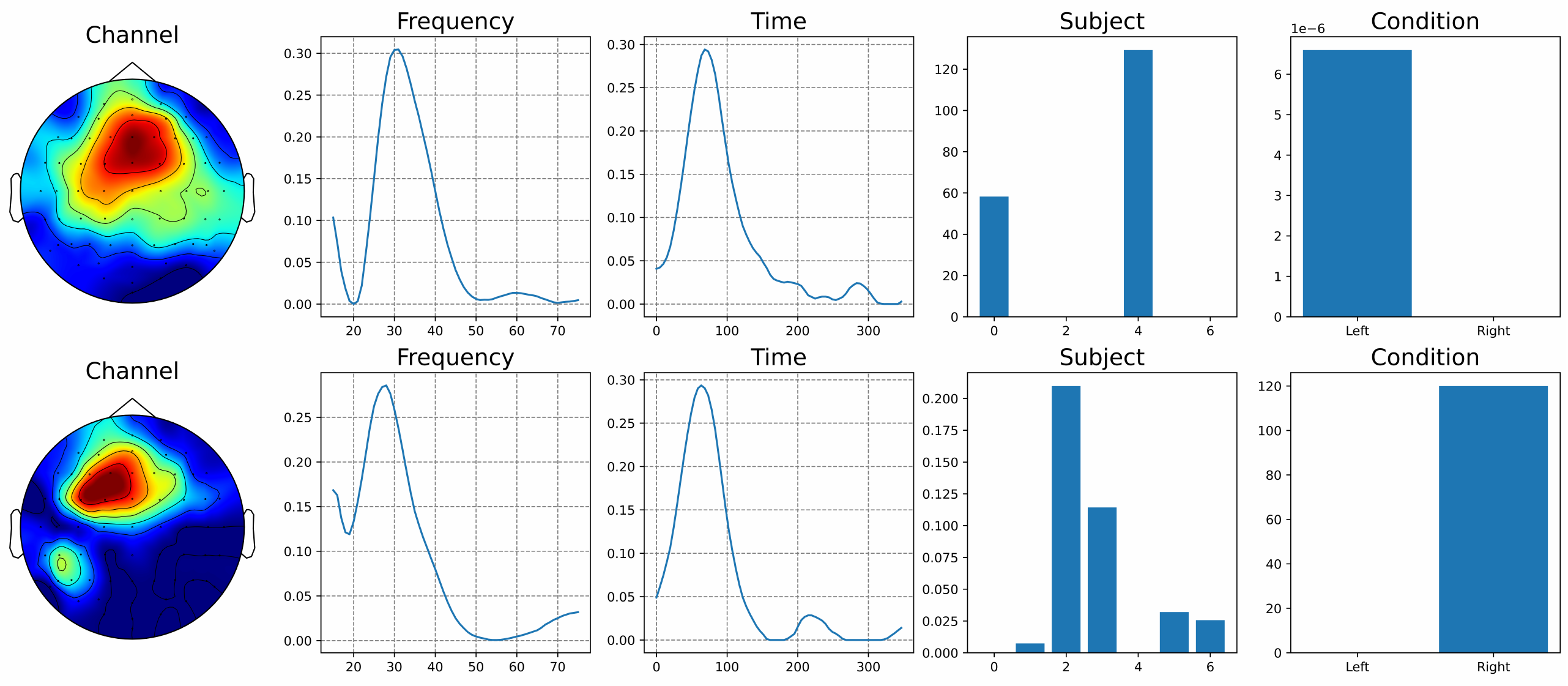}
    \label{fig5c}}
    \caption{Components of beta-band oscillations. (a) Activity appeared in the frontal lobe within 75 ms at 15 - 20 Hz. (b) activity occurs in the temporal lobe at 15 - 20 Hz, in the region of 75 ms. (c) Activity occurs in the frontal lobe at 25 - 30 Hz, in the region of 60 ms.}
    \label{fig5}
\end{figure*}

\begin{figure*}[htb]
    \centering
    \subfigtopskip=2pt
    \subfigbottomskip=2pt
    \subfigcapskip=-5pt
    \subfigure[]{\includegraphics[width=13cm]{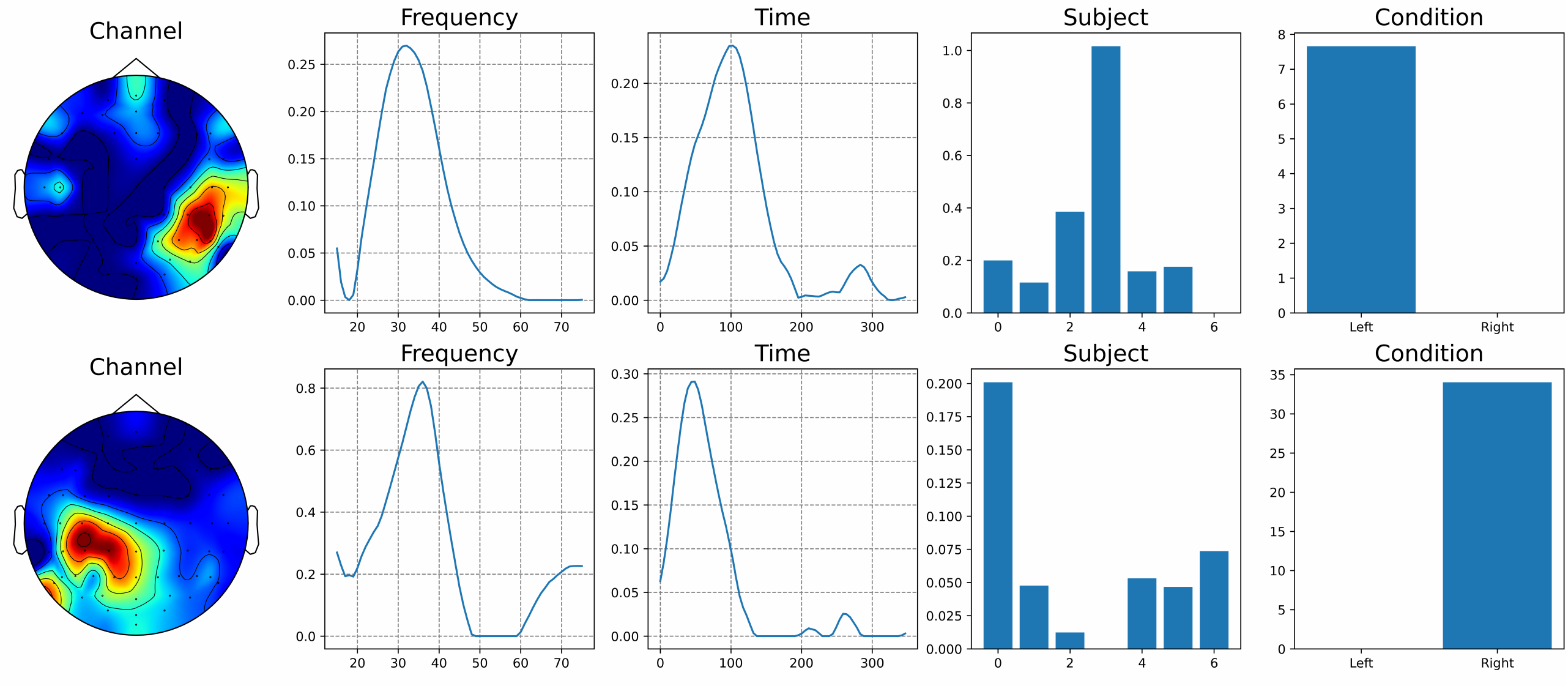}
    \label{fig6a}}
    \subfigure[]{\includegraphics[width=13cm]{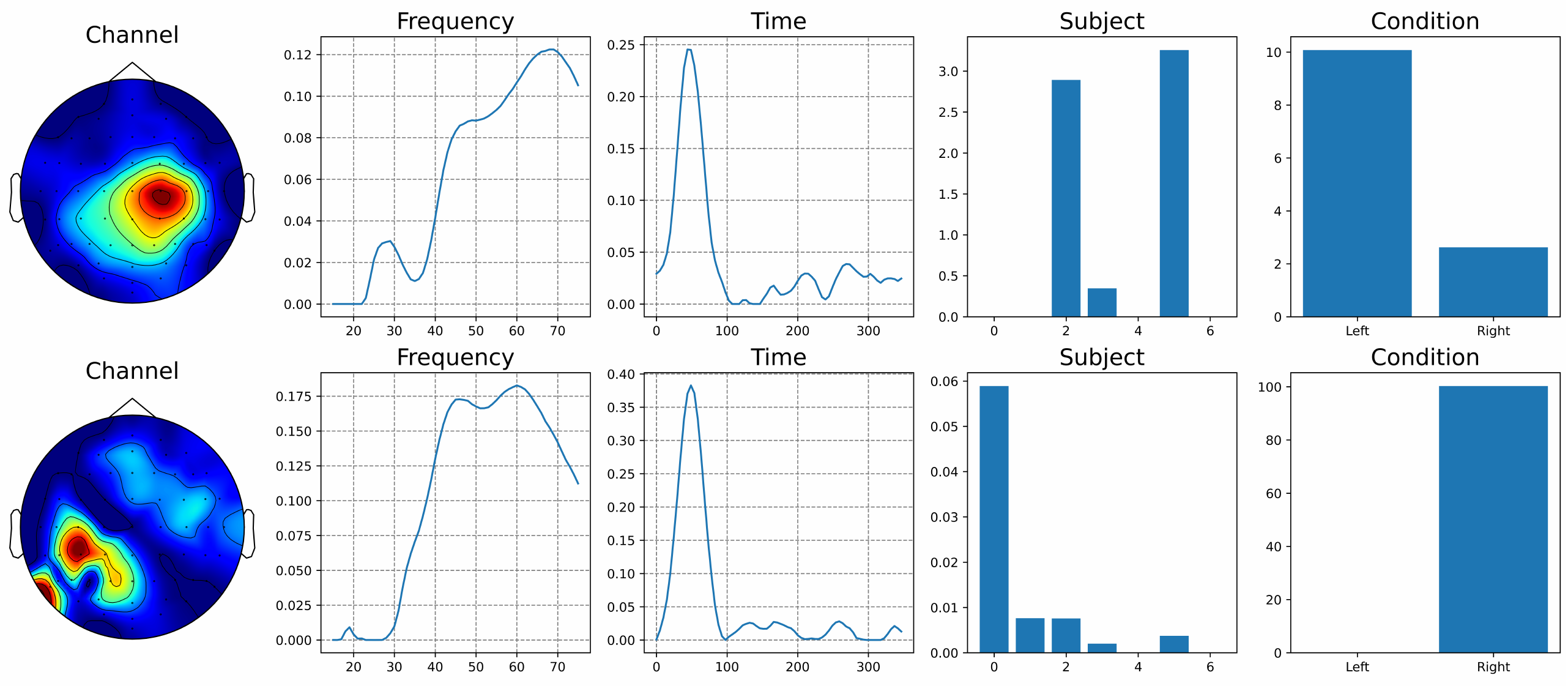}
    \label{fig6b}}
    \caption{Components of gamma-band oscillations. (a) Activity occurs within the frontal lobe at 30 - 40 Hz, In the region of  60 ms. (b) Activity in the first component occurs between the right parietal and temporal lobes at 40 - 75 Hz and 75 ms. Activity in the second component occurs between the left parietal and temporal lobes at 40 - 75 Hz, in the region of 45 ms.}
    \label{fig6}
\end{figure*}

\subsection*{Iteration stopping conditions}
Define the relative error/residual of FCNCP at the kth client, i-th iteration as:

\begin{equation}
    \text{RelErr}_{k}^{i}=\frac{{{\left\| {{\chi }_{k}}-\sum\nolimits_{r=1}^{{{R}_{k}}}{u_{r,k}^{\left( 1 \right),i}\circ u_{r,k}^{\left( 2 \right),i}\circ \cdots \circ u_{r.k}^{\left( N \right),i}} \right\|}_{F}}}{{{\left\| {{\chi }_{k}} \right\|}_{F}}}
\end{equation}

We terminate the FCNCP decomposition process when the following condition is satisfied between two iterations $i$ and $i + 1$.

\begin{equation}
    \left| \text{RelErr}_{k}^{i+1}-\text{RelErr}_{k}^{i} \right|<\varepsilon
\end{equation}\\

where $\varepsilon ={{10}^{-8}}$ is taken in this study.

\subsection*{Result analyse}
In this work, we compare with the conclusions drawn in the literature on non-negative tensor decomposition using all data\cite{2018Extracting}, and validate the effectiveness of the algorithm by comparing the consistency of the conclusions between the two. The results from the decomposition of the two clients were filtered for components of interest, and after integrating and visualising the results the following conclusions were drawn:

Figure~\ref{fig5} illustrates the components associated with $\beta$-band (15-30Hz) oscillations. In Figure~\ref{fig5a}, a pair of components displaying symmetrical activity in the left and right frontal lobes, triggered by contralateral hand stimulation, is depicted. Both components manifest within the 15-20 Hz range at approximately 70 ms. Remarkably, this identical pair of components was also observed in a previous study\cite{morup2006decomposing}. Figure~\ref{fig5b} showcases symmetrical activity in the left and right temporal lobes in the 15-20 Hz range in about 75 ms. The activity presented in Figure~\ref{fig5c} occurs in the 25-30 Hz range at around 60 ms within the frontal lobe. It is noteworthy that these $\beta$-band activities emerged within the first 100 ms following stimulus onset, aligning with earlier research on $\beta$ frequency oscillations in proprioceptive information processing\cite{2011Attenuation}.

Figure~\ref{fig6} introduces components that exhibit oscillations within the $\gamma$ band (30-75 Hz). Figure~\ref{fig6a} demonstrates symmetrical activity in the left and right temporal lobes at 30-40 Hz at approximately 60 ms. Interestingly, this pair of components was also identified in a previous study\cite{morup2006decomposing}. In Figure~\ref{fig6b}, two additional components with frequencies spanning 40-75 Hz are observed between the parietal and temporal lobes, appearing in about 75 ms and 45 ms, respectively. Notably, previous works\cite{morup2006decomposing,2011Attenuation} have argued that proprioceptive stimulation can induce gamma-band activity (GBA) within the 30-80 Hz range.

From the above analysis, it can be concluded that the coupled tensor decomposition using the federated learning idea to analyze the EEG signals can lead to the same conclusions as compared to the decomposition analysis of all the data. It can be proved that the method can effectively deal with high-order EEG data, and some key hidden information can be preserved.

\section*{Discussion}
To the best of our knowledge, this study represents the pioneering effort to analyze EEG signals using the coupled tensor decomposition federated learning method. In the context of current EEG signal processing methods, which are primarily tailored for data centralized on the same server, we conduct a deep consideration into the coupling relationships among data distributed across various servers. This relationship leads to the development of a novel coupled tensor decomposition model integrated within the federated learning framework. The model has good discriminative properties of tensor decomposition in high dimensional data representation and decomposition. Furthermore, it capitalizes on the unique advantages offered by coupled tensor decomposition for the cross-sample tensor data analysis. Simultaneously, it harnesses the strengths of federated learning for collaborative modeling on distributed servers. A significant achievement of this work is the resolution of the challenge posed by enabling joint analysis of data across servers without being able to share the data. We accomplish this by employing the elastic average federation framework to establish coupling constraints. This approach protects privacy and also enables individual clients to explore relevant information on their own while establishing coupling constraints across servers. Consequently, we offer a comprehensive suite of solutions for the efficient coupled tensor decomposition of high-dimensional EEG data distributed across different servers. Our methodology enables a profound exploration and extraction of similarities and distinctions within data collected from various subjects, providing new insights into EEG signal analysis in a distributed and privacy-preserving manner.

In our experiments, we initially utilized simulated tensor data measuring 61 × 72 × 64 (frequency, time, channel) while imposing coupling constraints on time and frequency. We then performed decomposition using the third-order FCNCP model. Our assessment of the stability and consistency of the decomposition results substantiates the algorithm's effectiveness. Subsequently, we applied a fifth-order FCNCP model to decompose a collection of fifth-order ERP tensor data. These data were acquired through the application of proprioceptive stimuli to the left and right hands. We executed the decomposition across all modalities (condition, subject, channel, frequency, time) while preserving the intermodal interaction information. Our experimental findings indicate that the conclusions drawn from the components decomposed by the FCNCP algorithm align with interpretations from relevant studies in cognitive neuroscience. This alignment underscores the method's capacity to effectively process higher-order EEG data while retaining essential hidden information.

However, we note that the current selection of coupling components based on the correlation coefficient matrix in the algorithm only applies to two clients, and the sparsity of the conclusions drawn for each client needs to be improved. Therefore, in future research, we will try to select coupling components for multiple clients using clustering. Meanwhile, we will also consider different optimisation methods (e.g. APG, MU, etc.) in the update of the client factor matrices and include sparsity constraints in the objective function to ensure the interpretability of the final results.

\section*{Conclusion}
Brain science currently holds a prominent position in international research, with many countries prioritizing its advancement. Collaborative efforts across servers prove instrumental in expediting progress within the field of brain science. In this study, we introduce a novel coupled tensor decomposition model, built upon the federated learning framework. This innovative approach leverages federated learning to establish coupling constraints for coupled tensor decomposition, effectively addressing the challenge of conducting joint data analysis across servers without being able to share the data. This accomplishment empowers the algorithm to delve deep into the similarities and differences among diverse subject data. Our experimental results demonstrate the successful application of this algorithm to both simulated data and ERP data collected through proprioceptive stimuli applied to the left and right hands. The outcomes affirm the algorithm's efficacy in processing higher-order EEG data while preserving crucial hidden information. This research not only equips us with new tools and methodologies for addressing the complexities of processing and analyzing high-dimensional EEG data on diverse servers but also advances the field of coupled tensor decomposition techniques and their integration with emerging federated learning frameworks. This contribution holds significant theoretical and practical value.

\section*{Acknowledgments}
This work was supported by Dalian Science and Technology Talent Innovation Support Project (No. 2023RY034).


\begin{thebibliography}{10}

\bibitem{1927The}
F.~L. Hitchcock.
\newblock The expression of a tensor or a polyadic as a sum of products.
\newblock {\em Journal of Mathematical Physics}, 6(1):164–189, 1927.

\bibitem{2016Tensor}
Nicholas~D Sidiropoulos, Lieven De~Lathauwer, Xiao Fu, Kejun Huang, Evangelos~E Papalexakis, and Christos Faloutsos.
\newblock Tensor decomposition for signal processing and machine learning.
\newblock 2016.

\bibitem{2015Tensor}
Andrzej Cichocki, Danilo Mandic, Lieven De~Lathauwer, Guoxu Zhou, Qibin Zhao, Cesar Caiafa, and Huy~Anh Phan.
\newblock Tensor decompositions for signal processing applications from two-way to multiway component analysis.
\newblock {\em IEEE Signal Processing Magazine}, 32(2):145--163, 2015.

\bibitem{2015Linked}
Guoxu Zhou, Qibin Zhao, Yu~Zhang, Tülay Adal, Shengli Xie, and Andrzej Cichocki.
\newblock Linked component analysis from matrices to high order tensors: Applications to biomedical data, 2015.

\bibitem{2011Applications}
Morten M{\o}Rup.
\newblock Applications of tensor (multiway array) factorizations and decompositions in data mining.
\newblock {\em Wiley Interdisciplinary Reviews: Data Mining and Knowledge Discovery}, 1(1):24--40, 2011.

\bibitem{1966Some}
Ledyard Tucker.
\newblock Some mathematical notes on three-mode factor analysis.
\newblock {\em Psychometrika}, 31(3):279--311, 1966.

\bibitem{Henk1999PARAFAC2}
Henk~A Kiers, Jos~M. Ten~Berge, and Rasmus Bro.
\newblock Parafac2—part i. a direct fitting algorithm for the parafac2 model.
\newblock {\em Journal of Chemometrics}, 12(3):223--238, 1999.

\bibitem{2020Tracing}
Marie Roald, Suchita Bhinge, Chunying Jia, Vince Calhoun, Tulay Adali, and Evrim Acar.
\newblock Tracing network evolution using the parafac2 model.
\newblock In {\em IEEE International Conference on Acoustics, Speech and Signal Processing}, 2020.

\bibitem{2009Nonnegative}
Andrzej Cichocki, Rafal Zdunek, Anh~Huy Phan, and Shun~Ichi Amari.
\newblock {\em Nonnegative Matrix and Tensor Factorizations: Applications to Exploratory Multi-Way Data Analysis and Blind Source Separation}.
\newblock Nonnegative Matrix and Tensor Factorizations: Applications to Exploratory Multi-Way Data Analysis and Blind Source Separation, 2009.

\bibitem{cichocki2007hierarchical}
Andrzej Cichocki, Rafal Zdunek, and Shun-ichi Amari.
\newblock Hierarchical als algorithms for nonnegative matrix and 3d tensor factorization.
\newblock In {\em International Conference on Independent Component Analysis and Signal Separation}, page 169–176. Springer, 2007.

\bibitem{cichocki2009fast}
Andrzej Cichocki and Anh~Huy Phan.
\newblock Fast local algorithms for large scale nonnegative matrix and tensor factorizations.
\newblock {\em IEICE transactions on fundamentals of electronics, communications and computer sciences}, 92(3):708–721, 2009.

\bibitem{boyd2011distributed}
Stephen Boyd, Neal Parikh, Eric Chu, Borja Peleato, and Jonathan Eckstein.
\newblock Distributed optimization and statistical learning via the alternating direction method of multipliers.
\newblock {\em Foundations and Trends in Machine learning}, 3(1):1–122, 2011.

\bibitem{kim2012fast}
Jingu Kim and Haesun Park.
\newblock Fast nonnegative tensor factorization with an active-set-like method.
\newblock In {\em High-Performance Scientific Computing}, page 311–326. Springer, 2012.

\bibitem{xu2015alternating}
Yangyang Xu.
\newblock Alternating proximal gradient method for sparse nonnegative tucker decomposition.
\newblock {\em Mathematical Programming Computation}, 7(1):39–70, 2015.

\bibitem{wang2021inexact}
Dongdong Wang and Fengyu Cong.
\newblock An inexact alternating proximal gradient algorithm for nonnegative cp tensor decomposition.
\newblock {\em Science China Technological Sciences}, 64(9):1893–1906, 2021.

\bibitem{2014Lowrank}
Fengyu Cong, Guoxu Zhou, Piia Astikainen, Qibin Zhao, Qian Wu, Asoke~K Nandi, Jari~K Hietanen, Tapani Ristaniemi, and Andrzej Cichocki.
\newblock Low-rank approximation based non-negative multi-way array decomposition on event-related potentials.
\newblock {\em International Journal of Neural Systems}, 24(8), 2014.

\bibitem{2016Linked}
Guoxu Zhou, Qibin Zhao, Yu~Zhang, Tulay Adali, Shengli Xie, and Andrzej Cichocki.
\newblock Linked component analysis from matrices to high-order tensors: Applications to biomedical data.
\newblock {\em Proceedings of the IEEE}, 104(2):310--331, 2016.

\bibitem{acar2013structure}
Evrim Acar, Andreas~J Lawaetz, Morten~A Rasmussen, and Rasmus Bro.
\newblock Structure-revealing data fusion model with applications in metabolomics. in embc (pp.6023-6026).
\newblock In {\em 35th Annual International Conference of the IEEE Engineering in Medicine and Biology Society (EMBC)}, pages 6023--6026, Osaka, Japan, July 2013.

\bibitem{2019Unraveling}
Evrim Acar, Carla Schenker, Yuri Levin-Schwartz, Vince~D. Calhoun, and Tülay Adali.
\newblock Unraveling diagnostic biomarkers of schizophrenia through structure-revealing fusion of multi-modal neuroimaging data.
\newblock {\em Frontiers in Neuroscience}, 13, 2019.

\bibitem{2020Group}
Xiulin Wang, Wenya Liu, Petri Toiviainen, Tapani Ristaniemi, and Fengyu Cong.
\newblock Group analysis of ongoing eeg data based on fast double-coupled nonnegative tensor decomposition.
\newblock {\em Journal of Neuroscience Methods}, 330:108502, 2020.

\bibitem{2022Exploring}
Wenya Liu, Xiulin Wang, Timo H$\ddot{a}$m$\ddot{a}$l$\ddot{a}$inen, and Fengyu Cong.
\newblock Exploring oscillatory dysconnectivity networks in major depression during resting state using coupled tensor decomposition.
\newblock {\em IEEE Transactions on Biomedical Engineering}, 69(8):2691--2700, 2022.

\bibitem{hunyadi2017tensor}
Borbála Hunyadi, Patrick Dupont, Wim Van~Paesschen, and Sabine Van~Huffel.
\newblock Tensor decompositions and data fusion in epileptic electroencephalography and functional magnetic resonance imaging data.
\newblock {\em Wiley Interdisciplinary Reviews: Data Mining and Knowledge Discovery}, 7(1):e1197, 2017.

\bibitem{martinez2004concurrent}
Eduardo Mart{'\i}nez-Montes, Pedro~A Vald{'e}s-Sosa, Fumikazu Miwakeichi, Robin~I Goldman, and Mark~S Cohen.
\newblock Concurrent eeg/fmri analysis by multiway partial least squares.
\newblock {\em Neuroimage}, 22(3):1023–1034, 2004.

\bibitem{acar2017tensor}
Evrim Acar, Yuri Levin-Schwartz, Vince~D Calhoun, and T{"u}lay Adali.
\newblock Tensor-based fusion of eeg and fmri to understand neurological changes in schizophrenia.
\newblock In {\em 2017 IEEE International Symposium on Circuits and Systems (ISCAS)}, page 1–4. IEEE, 2017.

\bibitem{jonmohamadi2020extraction}
Yaqub Jonmohamadi, Suresh Muthukumaraswamy, Jian Chen, James Roberts, Ross Crawford, and Avinash Pandey.
\newblock Extraction of common task features in eeg-fmri data using coupled tensor-tensor decomposition.
\newblock {\em Brain Topography}, 33(5):636–650, 2020.

\bibitem{shi2016edge}
Weisong Shi, Jie Cao, Quan Zhang, Youhuizi Li, and Lanyu Xu.
\newblock Edge computing: Vision and challenges.
\newblock {\em IEEE internet of things journal}, 3(5):637–646, 2016.

\bibitem{yu2020new}
Hyeonseok Yu and Young-Ho Kim.
\newblock New rsa encryption mechanism using one-time encryption keys and unpredictable bio-signal for wireless communication devices.
\newblock {\em Electronics}, 9(2):246, 2020.

\bibitem{li2020review}
Liangzhen Li, Yifan Fan, Michael Tse, and Kaixuan~Yvonne Lin.
\newblock A review of applications in federated learning.
\newblock {\em Computers \& Industrial Engineering}, 149:106854, 2020.

\bibitem{li2020federated}
Tian Li, Anit~Kumar Sahu, Ameet Talwalkar, and Virginia Smith.
\newblock Federated learning: Challenges, methods, and future directions.
\newblock {\em IEEE Signal Processing Magazine}, 37(3):50–60, 2020.

\bibitem{zhang2021survey}
Chaoyang Zhang, Yuxiang Xie, Haoxiang Bai, Bo~Yu, Wenlin Li, and Yang Gao.
\newblock A survey on federated learning.
\newblock {\em Knowledge-Based Systems}, 216(1):106775, 2021.

\bibitem{mammen2021federated}
Mammen P.~Mammen Jr.
\newblock Federated learning: Opportunities and challenges.
\newblock 2021.

\bibitem{yang2019federated}
Qiang Yang, Yang Liu, Tianjian Chen, and Yongxin Tong.
\newblock Federated machine learning: Concept and applications.
\newblock {\em ACM Transactions on Intelligent Systems and Technology (TIST)}, 10(2):1–19, 2019.

\bibitem{abdulrahman2020survey}
Sanaa AbdulRahman, Hala Tout, Hamza Ould-Slimane, Azzam Mourad, Chamseddine Talhi, and Mohsen Guizani.
\newblock A survey on federated learning: The journey from centralized to distributed on-site learning and beyond.
\newblock {\em IEEE Internet of Things Journal}, 8(7):5476–5497, 2020.

\bibitem{mcmahan2017communication}
H~Brendan McMahan, Eider Moore, Daniel Ramage, Seth Hampson, and Blaise~Aguera y~Arcas.
\newblock Communication-efficient learning of deep networks from decentralized data.
\newblock In {\em Artificial intelligence and statistics}, page 1273–1282. PMLR, 2017.

\bibitem{kong2019federated}
Linghe Kong, Xiao-Yang Liu, Huanhuan Sheng, Peng Zeng, and Guihai Chen.
\newblock Federated tensor mining for secure industrial internet of things.
\newblock {\em IEEE Transactions on Industrial Informatics}, 16(3):2144–2153, 2019.

\bibitem{kim2017federated}
Youngduck Kim, Jimeng Sun, Hwanjo Yu, and Xiaoqian Jiang.
\newblock Federated tensor factorization for computational phenotyping.
\newblock In {\em Proceedings of the 23rd ACM SIGKDD International Conference on Knowledge Discovery and Data Mining}, page 887–895, 2017.

\bibitem{gao2021federated}
Yuxiang Gao, Guoqi Zhang, Chaoyang Zhang, Jianqiang Wang, Laurence~T. Yang, and Yong Zhao.
\newblock Federated tensor decomposition-based feature extraction approach for industrial iot.
\newblock {\em IEEE Transactions on Industrial Informatics}, 17(12):8541–8549, 2021.

\bibitem{zhang2015deep}
Sixin Zhang, Anna~E Choromanska, and Yann LeCun.
\newblock Deep learning with elastic averaging sgd.
\newblock {\em Advances in neural information processing systems}, 28, 2015.

\bibitem{morup2007erpwavelab}
Morten M{\o}rup, Lars~Kai Hansen, and Sidse~M Arnfred.
\newblock Erpwavelab: a toolbox for multi-channel analysis of time–frequency transformed event related potentials.
\newblock {\em Journal of neuroscience methods}, 161(2):361–368, 2007.

\bibitem{morup2006decomposing}
Morten M{\o}rup, Lars~Kai Hansen, Josef Parnas, and Sidse~M Arnfred.
\newblock Decomposing the time-frequency representation of eeg using non-negative matrix and multi-way factorization.
\newblock Technical report, Technical University of Denmark Technical Report, 2006.

\bibitem{delorme2004eeglab}
Arnaud Delorme and Scott Makeig.
\newblock Eeglab: an open source toolbox for analysis of single-trial eeg dynamics including independent component analysis.
\newblock {\em Journal of neuroscience methods}, 134(1):9–21, 2004.

\bibitem{david2006mechanisms}
Olivier David, James~M Kilner, and Karl~J Friston.
\newblock Mechanisms of evoked and induced responses in meg/eeg.
\newblock {\em Neuroimage}, 31(4):1580–1591, 2006.

\bibitem{cohen2014analyzing}
Mike~X Cohen.
\newblock {\em Analyzing neural time series data: theory and practice}.
\newblock MIT press, 2014.

\bibitem{2018Extracting}
Wang Deqing, Zhu Yongjie, Ristaniemi Tapani, and Cong Fengyu.
\newblock Extracting multi-mode erp features using fifth-order nonnegative tensor decomposition.
\newblock {\em Journal of Neuroscience Methods}, pages 240--247, 2018.

\bibitem{2011Attenuation}
Sidse Marie~Hemmingsen Arnfred, Morten M{\o}Rup, J{\o}Rgen Thalbitzer, Lennart Jansson, and Josef Parnas.
\newblock Attenuation of beta and gamma oscillations in schizophrenia spectrum patients following hand posture perturbation.
\newblock {\em Psychiatry Research}, 185(1-2):215--224, 2011.

\end{thebibliography}

\bibliographystyle{abbrv}

\end{document}